\documentclass[english,12pt]{amsart}
\usepackage[utf8x]{inputenc}
\usepackage{ucs}
\usepackage[top=1in, bottom=1in, left=1in, right=1in]{geometry}
\usepackage{amsmath}
\usepackage{amsthm}
\usepackage{stmaryrd}
\usepackage{amsfonts,dsfont}
\usepackage{changes}
\usepackage{amssymb}
\usepackage{mathtools} 
\usepackage{amsthm}
\usepackage{framed}
\usepackage{mathrsfs} 
\usepackage{graphicx}
\usepackage{enumerate}  

\newtheorem{theoA}{Theorem}

\newtheorem{Theorem}{Theorem}

\newtheorem{Corollary}[Theorem]{Corollary}

\newtheorem{Remark}[Theorem]{Remark}

\graphicspath{ {./images/} }  
\renewcommand{\sinh}{\operatorname{sh}}
\renewcommand{\cosh}{\operatorname{ch}}

\title{On \(H^{2|2}\) isomorphism theorems and reinforced loop soup}
\author[Y. Chang]{Yinshan Chang}
\address{College of Mathematics, Sichuan University, Chengdu 610065, China}
\email{ychang@scu.edu.cn}
\author[DZ. Liu]{Dang-Zheng Liu}
\address{Key Laboratory of Wu Wen-Tsun Mathematics, CAS, School of Mathematical Sciences, University of Science and Technology of China, Hefei 230026, P.R.~China}
\email{dzliu@ustc.edu.cn}
\author[X. Zeng]{Xiaolin Zeng}
\address{Institut de Recherche Mathématique Avancée, UMR 7501, Université de Strasbourg, 7 rue René Descartes, 67084 Strasbourg, France}
\email{zeng@math.unistra.fr}
\thanks{The first author is supported by NSFC Grant 11701395, the second author is supported by National Natural Science Foundation of China \#11771417, the Youth Innovation Promotion Association CAS \#2017491.}
\begin{document}
\maketitle
\begin{abstract}
We show that supersymmetric (susy) hyperbolic isomorphism theorems that relate Vertex Reinforced Jump Processes and \(H^{2|2}\) field, introduced in \cite{bauerschmidt2018dynkin} and \cite{bauerschmidt2019geometry}, are annealed version of isomorphism theorems relating Markov processes and Gaussian free field, with the help of a Bayes formula that relates susy hyperbolic field to susy free field. On the other hand, we also prove a BFS--Dynkin's isomorphism theorem for reinforced loop soup. Moreover, we provide yet another proof of BFS-Dynkin's isomorphism for VRJP à la Feynman--Kac.
\end{abstract}
\section{Introduction}
\label{sec:org6c87435}

In \cite{bauerschmidt2018dynkin} and \cite{bauerschmidt2019geometry}, several isomorphism theorems relating a susy hyperbolic sigma field (called the \(H^{2|2}\) model, see \cite{DSZ06} for details) and the Vertex Reinforced Jump Process (VRJP, discussed in e.g. \cite{basdevant2012continuous,davis2002continuous,davis2004vertex,STZ15,SZ15,sabot2019polynomial,merkl2019random,poudevigne2019monotonicity,collevecchio2018note} etc list non exclusive, this process is also related to the edge reinforced random walk, e.g. see references in \cite{SZ15}) are introduced, and were applied to deduce the recurrence of the VRJP on two dimensional lattice. In this paper we try to answer the following question: \textbf{\textbf{are these isomorphism theorems related to the classical ones (e.g. \cite{eisenbaum1994dynkin},\cite{sabot2016inverting},\cite{le2008dynkin}) that correspond to Gaussian free field?}}

We provide a positive answer to this question and give alternative proofs to these isomorphism theorems. In particular, we prove a Bayes formula (Theorem \ref{orgc9bde69}) that relates the susy hyperbolic sigma field to the classical Gaussian free field. Using the Bayes formula, we also provide the susy hyperbolic version of BFS--Dynkin's isomorphism for loop soup (Theorem \ref{org6c0a237}), which we define as reinforced loop soup. For motivation on loop soup, see \cite{symanzik1968euclidean} and \cite{le2010markov} and references therein. Finally, we also give another proof of susy hyperbolic BFS--Dynkin's isomorphism à la Feynman-Kac in the appendix.
\subsection{Definitions and notations}
\label{sec:orgf8d5325}
In this section, we define our notations and gather some useful facts about our subjects. We will exclusively work on a finite graph, denoted by \(\mathcal{G}=(V,E)\), where the edges are non-oriented, to each edge \(\{i,j\}\in E\) we associate a positive real \(W_{i,j}\), and \(W_{i,j}=0\) if there is no edge between \(i,j\). This weighted graph is encoded by the following matrix
\[\Delta_W=(\Delta_W(i,j))_{i,j\in V},\ \text{ where }\Delta_W(i,j)=\begin{cases}-W_{i,j} & i\ne j \\ \sum_{k\in V\setminus\{i\}}W_{i,k} & i=j\end{cases}.\]
Sometimes we add a cemetery point to our graph, denoted by \(\delta\), so the graph with the cemetery is \(\widetilde{G}=(\widetilde{V},\widetilde{E})\), where \(\widetilde{V}=V\cup \{\delta\}\) and \(\widetilde{E}=E\cup \{\{\delta,i\}: W_{\delta,i}>0\}\). We assume that at least one of the \(W_{\delta,i}>0\) (so the graph is connected), usually all the \(W_{\delta,i}\) are the same, in such case, we usually denote \(W_{\delta,i}=h,\ \forall i\in V\). This enlarged graph is encoded by
\[\Delta_{\widetilde{W}}=(\Delta_{\widetilde{W}}(i,j))_{i,j\in \widetilde{V}},\ \text{ where }\Delta_{\widetilde{W}}(i,j)=\begin{cases}-W_{i,j} & i\ne j \\ \sum_{k\in \widetilde{V}\setminus\{i\}}W_{i,k} & i=j\end{cases}.\]
Note that we have \((\Delta_{\widetilde{W}})_{V\times V}=\Delta_W+h\), where \((\Delta_{\widetilde{W}})_{V\times V}\) is the restriction of \(\Delta_{\widetilde{W}}\) to the set \(V\) and \(h\) is considered as the diagonal matrix \(h \operatorname{Id}\).
\subsubsection{The susy hyperbolic sigma model}
\label{sec:org75ab0ab}
The \(H^{2|2}\) model is discussed in detail in \cite{DSZ06}, in this subsection we briefly recall the basic properties of this model. {For more details, we refer to \cite[Section~2]{bauerschmidt2018dynkin} and the reference therein.} To each vertex \(i\in V\), we associate a (super) vector \(\Phi_i=(x_i,y_i,z_i,\xi_i,\eta_i)\). {Here, \(x_i\) and \(y_i\) are real variables, \(\xi_i\) and \(\eta_i\) are Grassmann variables. Considerable care must be taken when dealing with Grassmann variables, the product between \(\xi_i\) and \(\eta_j\) are anticommutative, i.e. \(\xi_i\eta_j=-\eta_j\xi_i\). In particular, \(\xi_i\xi_i=\xi_i^2=0\).} This vector \(\Phi_i\) lives in the \(H^{2|2}\) susy hyperbolic space, i.e.
\[\Phi_i\cdot \Phi_i=\Phi_i^2=x_i^2+y_i^2-z_i^2+2\xi_i\eta_i=-1,\]
where \(z_i\) is chosen in the positive branch, i.e.
\[z_i=\sqrt{1+x_i^2+y_i^2+2\xi_i\eta_i}{=\sqrt{1+x_i^2+y_i^2}+\frac{\xi_i\eta_i}{\sqrt{1+x_i^2+y_i^2}}}.\]

{Here, one should think of \(z_i=\sqrt{1+x_i^2+y_i^2+2\xi_i\eta_i}\) as a super-function of \(x_i\), \(y_i\), \(\xi_i\) and \(\eta_i\). Its precise meaning is given by the Taylor expansion of \(\sqrt{\cdot}\) at \(1+x_i^2+y_i^2\) as follows: \[z_i=\sqrt{1+x_i^2+y_i^2}+\frac{\xi_i\eta_i}{\sqrt{1+x_i^2+y_i^2}}-\frac{\xi_i\eta_i\xi_i\eta_i}{2(1+x_i^2+y_i^2)^{\frac{3}{2}}}+\cdots.\]

Since \(\xi_i\xi_i=0\), we have \(\xi_i\eta_i\xi_i\eta_i=-\xi_i\xi_i\eta_i\eta_i=0\), therefore the third term vanishes. Similarly, all the higher order terms vanish and the expansion has only the first two terms.}

The collection of vectors \((\Phi_i)_{i\in V}=\Phi\) live in \((H^{2|2})^V\), on which the inner product of two vectors is defined by
\[\Phi_i\cdot \Phi_j=\Phi_i\Phi_j =x_ix_j+y_iy_j-z_iz_j+\xi_i\eta_j-\eta_i\xi_j.\]
This inner product is symmetric, i.e. \(\Phi_i\Phi_j = \Phi_j\Phi_i\), the zero vector is \((0,0,1,0,0)\), which is simply denoted by \(0\). The energy of the \(H^{2|2}\) model on \(V\) is defined by
\begin{equation}
\label{eq-h22-energy}
\begin{aligned}
\frac{1}{2}\Phi \Delta_W\Phi&=\frac{1}{2}\Phi^t \Delta_W\Phi=\frac{1}{2}\sum_{i,j\in V}(\Delta_W)_{ij}\Phi_i\cdot\Phi_j\\
&=-\sum_{\{i,j\}\in E} W_{i,j}(x_ix_j+y_iy_j-z_iz_j+\xi_i\eta_j+\xi_j\eta_i+1).
\end{aligned}
\end{equation}
The Haar measure (Berezin integral form, rather) on \((H^{2|2})^V\) equals
\begin{equation}
\label{eq-h22-haar}
D\mu_V(\Phi)=\prod_{i\in V} \frac{1}{2\pi z_i}dx_i dy_i d\xi_i d\eta_i.
\end{equation}
{
Later, we will integrate a super-function \(F(x,y,\xi,\eta)\) against \(D\mu_V\) and denote \(\int F D\mu_{V}(\Phi)\) the resulting integral. Equivalently, this amounts to integrate the super-function \(G=\frac{F}{\prod_{i\in V}z_i}\) against \(\prod_{i\in V}\frac{dx_idy_id\xi_id\eta_i}{2\pi}\). For the integration of a super-function \(F\) against \(\prod_{i\in V}\frac{dx_idy_id\xi_id\eta_i}{2\pi}\), we refer to \cite[Eqs. (28) and (29)]{bauerschmidt2018dynkin}.
}

Due to the non compact feature of hyperbolic space, the integral \(\int e^{-\frac{1}{2}  \Phi \Delta_W\Phi }D\mu_V(\Phi)\) is not normalizable, we have to add a pinning, or equivalently, a boundary condition. {i.e.}, we consider the graph \(\widetilde{V}\) instead, if at least one of the \(W_{\delta,i}>0\), then the following integral is normalizable:
\[\int \mathds{1}_{\Phi_{\delta}=0} e^{-\frac{1}{2}  \widetilde{\Phi}\Delta_{\widetilde{W}}\widetilde{\Phi} }D\mu_{V}(\Phi)=1,\]
where \(\widetilde{\Phi}=(\Phi_i)_{i\in \widetilde{V}}\). Note that
{
\begin{equation}
\label{eq-h22-energy-mass}
\widetilde{\Phi}\Delta_{\widetilde{W}}\widetilde{\Phi}=\sum_{i,j\in\widetilde{V}}(\Delta_{\widetilde{W}})_{ij}\widetilde{\Phi}_i\widetilde{\Phi}_{j}=-2\sum_{\{i,j\}\in \widetilde{E}} W_{i,j} (x_ix_j+y_iy_j-z_iz_j+\xi_i\eta_j+\xi_j\eta_i+1).
\end{equation}
In particular, if we have \(\Phi_{\delta}=0=(0,0,1,0,0)\), then
\begin{equation}
\widetilde{\Phi}\Delta_{\widetilde{W}}\widetilde{\Phi}=-2 \sum_{\{i,j\}\in E}W_{i,j}(x_ix_j+y_iy_j-z_iz_j+\xi_i\eta_j+\xi_j\eta_i+1)-2\sum_{i\in V}W_{\delta,i}(-z_i+1)
\end{equation}
}

We recover the usual definition of \(H^{2|2}\)-model defined in \cite{DSZ06} when \(W_{\delta,i}=h\) for some \(h>0\) for all \(i\in V\), and in the language of statistical mechanics, \(h\) is the mass of the model. In the sequel, this is our default choice if we do not mention anything else explicitly, and is called \(H^{2|2}\) model on \(\widetilde{V}\) pinned at \(\Phi_{\delta}=0\), or equivalently, \(H^{2|2}\) model on \(V\) with mass \(h\). We denote
\begin{equation}
\label{eq-defn-tilde-W-with-zero-boundary}
\left< \cdot \right>_{\widetilde{W},\Phi_{\delta}=0}=\int \cdot \mathds{1}_{\Phi_{\delta}=0} e^{-\frac{1}{2} \widetilde{\Phi}\Delta_{\widetilde{W}}\widetilde{\Phi}} D\mu_V(\Phi)
\end{equation}
the expectation w.r.t. the above density measure, moreover, abusively, we denote
\begin{equation}
\label{eq-defn-W-without-boundary}
\left< \cdot \right>_W=\int \cdot e^{-\frac{1}{2} \Phi \Delta_W\Phi}D\mu_V(\Phi)
\end{equation}
the (non-normalizable) expectation, for which we must be very careful on the observable we integrate with.

There are some natural symmetries on this model, in particular, we will use two of them, which we recall here. The first one is called the supersymmetry (or \(Q\)-symmetry), which states that, if \(F=F(\Phi)\) is any function that is invariant under the symmetric group of the inner product on \(H^{2|2}\) space, then
\begin{equation}
\label{eq-h22-localized}
\left< F \right>_{\widetilde{W},\Phi_{\delta}=0}=F(0),
\end{equation}
where \(F(0)\) is simply the value obtained by taking its argument \(\Phi_i=0\) for all \(i\in V\). In particular, if \(F=F(z)\) is a function of the component \(z\), then \(\left< F \right>_{\widetilde{W},\Phi_{\delta}=0}=F(1)\).

The second type of invariance we want to discuss here is the rotational symmetric in the \(xy\) plane and Lorentz boost in the \(xz\) plane, which are defined (respectively) by
\begin{enumerate}
\item Euclidean rotation \(R_{\alpha}\), with \(\alpha\in [0,2\pi]\)
\[R_{\alpha}\Phi=(x\cos \alpha +y\sin \alpha,-x \sin \alpha+y\cos \alpha,z,\xi,\eta)\]
\item Lorentz boost \(\theta_s\), with \(s\in \mathbb{R}\),
\[\theta_s\Phi=(x\cosh s+z\sinh s , y,z\cosh s+x\sinh s,\xi,\eta).\]
\end{enumerate}
Both the rotation and the boost leave the inner product \(\left< \Phi_i,\Phi_j \right>\) and the Haar measure \(D\mu_V(\Phi)\) invariant, as a consequence, \(\left< \cdot \right>_W\) is also invariant. Note that if we boost with \(\left< \cdot \right>_{\widetilde{W},\Phi_{\delta}=0}\), the boundary condition is shifted, so \(\left< \cdot \right>_{\widetilde{W},\Phi_{\delta}=0}\) is not invariant. To be more precise, if we boost with \(\left< \cdot \right>_{\widetilde{W},\Phi_{\delta}=0}\), we obtain \(\left< \cdot \right>_{\widetilde{W},\Phi_{\delta}=\theta_s0}\), where { \(\theta_s 0=(\sinh s,0,\cosh s,0,0)\) and }

\begin{equation}
\label{eq-defn-tilde-W-with-s-boundary}
\left< \cdot  \right>_{\widetilde{W},\Phi_{\delta}=\theta_s 0}=\int \cdot \mathds{1}_{\Phi_{\delta}=\theta_s 0}e^{-\frac{1}{2} \widetilde{\Phi}\Delta_{\widetilde{W}} \widetilde{\Phi}} D\mu_V(\Phi).
\end{equation}
Similarly, for \(a\in V\) and \(s\in \mathbb{R}\), we define
\begin{equation}
\label{eq-defn-W-with-s-boundary}
\left< \cdot \right>_{W,\Phi_a=\theta_s 0}=\int \cdot \mathds{1}_{\Phi_a=\theta_s 0} e^{-\frac{1}{2} \Phi \Delta_W \Phi} D\mu_{V\setminus\{a\}}(\Phi).
\end{equation}
\subsubsection{Vertex Reinforced Jump Process}
\label{sec:orgfb60612}
The Vertex Reinforced Jump Process (VRJP) is first discussed in \cite{davis2002continuous}, most of the facts we recall here can be found in \cite{STZ15},\cite{ST15}. The VRJP on \(\mathcal{G}\) with edge weight \(W\) and initial local time \(z\in \mathbb{R}_{>0}^V\) is a continuous time jump process, denoted by \(Y=(Y_t)_{t\ge 0}\), starts from a vertex \(i_0\in V\), and it jumps from \(i\) to \(j\) at time \(t\) at rate \(W_{i,j}L_j(t)\), where
\[L_j(t):= z_j+\int_0^t \mathds{1}_{Y_s=j}ds .\]
The quantity \(L_j(t)\) is called its local time. If we perform the following time scaling: let
\[D(t)=\sum_{i\in V}\left( L_i(t)^2-z_i^2 \right)\]
and define \(Z_s=Y_{D^{-1}(s)}=Y_t\), where \(D(t)=s\), then we have
\begin{equation}
\label{eq-time-S-L}
S_i(s)=\int_0^s \mathds{1}_{Z_t=i}dt=L_i(D^{-1}(s))^2-z_i^2.
\end{equation}
The jump rate from \(i\) to \(j\) of the process \(Z=(Z_s)_{s\ge 0}\) at time \(s\) equals
\[\frac{1}{2} W_{i,j} \sqrt{\frac{S_j(s)+z_j^2}{S_i(s)+z_i^2}}.\]
The process \(Z\) turns out to be a mixture of Markov jump process, that is, to sample \(Z\), we can first sample the environment \(u\in\{u_i\in \mathbb{R},i\in V,\ u_{i_0}=0\}\), according to the probability density function
\begin{equation}
\label{eq-nu-W-z}
\nu_{i_0}^{W,z}(du)=\mathds{1}_{u_{i_0}=0}e^{-\frac{1}{2} \sum_{\{i,j\}\in E} W_{i,j}(e^{u_i-u_j}z_j^2+e^{u_j-u_i}z_i^2-2z_iz_j)}\sqrt{D(W,u)}\prod_{i\ne i_0}\frac{z_i e^{-u_i}}{\sqrt{2\pi}}du_i,
\end{equation}
where
\[D(W,u)=\sum_{\text{ spanning tree }T \text{ of }\mathcal{G}}\prod_{\{i,j\}\in T}W_{i,j}e^{u_i+u_j}.\]
Then, sample a Markov jump process (called the quenched process) with (static) jump rate from \(i\) to \(j\):
\[\frac{1}{2} W_{i,j}e^{u_j-u_i},\]
We denote \(E_{a}^u\) the expectation w.r.t. the probability law of the quenched process \((Z_s)\) starting from \(a\). The generator of the quenched process is not symmetric, but it can easily be made symmetric by a simple time change. First, let \(A^u=(A^u_{i,j})_{i,j\in V}\) be the matrix with entries
\begin{equation}\label{eq-defn-Au}
 A^u_{i,j}=\begin{cases}-W_{i,j} & i\ne j \\ \sum_{k\in V\setminus\{i\}}W_{i,k}e^{u_k-u_i} & i=j\end{cases}
\end{equation}
The generator of the quenched process equals \(\frac{1}{2} e^{-u}A^u e^u\), where \(e^u\) is considered as the \(V\times V\) diagonal matrix with diagonal entry \(e^{u_i}\). If we scale the time at each vertex by a factor \(2e^{2u_i}\), then we get a reversible Markov process \(\widehat{Z}_t\) with jump rate \(W_{i,j}e^{u_i+u_j}\), its generator is denoted by \(B^u\), where
\[B^u_{i,j}=\begin{cases} -W_{i,j}e^{u_i+u_j} & i\ne j \\ \sum_{k\in V\setminus\{i\}}W_{i,k}e^{u_i+u_k} & i=j.\end{cases}\]
The local time \(\widehat{S}_i(t)\) of \(\widehat{Z}_t\) is related to the local time \(S(t)\) of \(Z_t\) by \(2e^{2u_i}\widehat{S}_i(t)=S_i(t)\).

Similarly, if we consider the VRJP on \(\widetilde{V}\), let \(\zeta=\inf\{t>0:\ Y_t=\delta\}\), we kill the process when it hits \(\delta\). The process without killing is a mixture of Markov process for the same reason as above, the mixing measure is
\begin{equation}
\label{eq-nutilde-W-z}
\nu_{i_0}^{\widetilde{W},z}(du)=\mathds{1}_{u_{i_0}=0}e^{-\frac{1}{2}\sum_{\{i,j\}\in \widetilde{E}} W_{i,j}(e^{u_i-u_j}z_j^2+e^{u_j-u_i}z_i^2-2z_iz_j)}\sqrt{D(W,u)}\prod_{i\in \widetilde{V}\setminus \{i_0\}}\frac{z_i e^{-u_i}}{\sqrt{2\pi}}du_i,
\end{equation}
and we can also define \(\widetilde{A}^u\) and \(\widetilde{B}^u\) on \(\widetilde{V}\) as before for suitable time changes:
{
\begin{equation}\label{eq-defn-tilde-Au-and-tilde-Bu}
 \widetilde{A}^u_{i,j}=\begin{cases}-W_{i,j} & i\ne j\in\widetilde{V} \\ \sum_{k\in \widetilde{V}\setminus\{i\}}W_{i,k}e^{u_k-u_i} & i=j\in\widetilde{V},\end{cases}\,\quad \widetilde{B}^u_{i,j}=\begin{cases} -W_{i,j}e^{u_i+u_j} & i\ne j\in\widetilde{V} \\ \sum_{k\in \widetilde{V}\setminus\{i\}}W_{i,k}e^{u_i+u_k} & i=j\in\widetilde{V}.\end{cases}
\end{equation}
}
As a consequence, the killed VRJP is a mixture of Markov process, killed (in the quenched meaning) at its first hitting of \(\delta\).

Another property that we will use in our argument is the following. If we perform the change of variable \(v_i=u_i-u_b\) for \(i\in V\) in the measure \(\nu_a^{W,1}(du)\), then we have
\begin{equation}
\label{eq-atob}
e^{u_b-u_a}\nu_a^{W,1}(du)=\nu_b^{W,1}(dv).
\end{equation}
Hence, if \(F\) is a function of the gradients \((u_i-u_j:\ i,j\in V)\), then
\begin{equation}
\label{eq-change-root}
\int F e^{u_b-u_a} \nu_a^{W,1}(du)=\int F \nu_b^{W,1}(du).
\end{equation}
\subsubsection{Supersymmetric free field and Parisi-Sourlas formula}
\label{sec:orgb6dd39a}
In this subsection we discuss the classical supersymmetric free field, and its related Parisi-Sourlas formula (or localization of susy integral), and we apply the formula to show that \(H^{2|2}\) model shares the very same localization formula. This is a classical topic in supersymmetry, cf.~\cite{efetov1999supersymmetry},\cite{wegner2016supermathematics} etc.

The susy free field on \(V\) with generator \(W\) is a random vector \(\Psi=(\Psi_i)_{i\in V}\), where \(\Psi_i=(x_i,y_i,\xi_i,\eta_i)\) is a four-vector, with two real component and two fermions. Its energy is defined by
\begin{equation}
\label{eq-freefield-energy}
\frac{1}{2}  \Psi \Delta_W \Psi =\frac{1}{2}\sum_{i,j\in V}(\Delta_W)_{i,j}\Psi_i\Psi_j=\sum_{\{i,j\}\in E}W_{i,j}(\Psi_i-\Psi_j)(\Psi_i-\Psi_j)
\end{equation}
where
\[\Psi_i\Psi_j =x_ix_j+y_iy_j+\xi_i\eta_j+\xi_j\eta_i.\]
The Haar measure \(D\Psi\) is defined by
\begin{equation}
\label{eq-freefield-haar}
D\Psi=\prod_{i\in V}\frac{1}{2\pi}dx_i dy_i d\xi_i d\eta_i.
\end{equation}
{For the integration of a super-function against \(D\Psi\), we refer to \cite[Eqs. (28) and (29)]{bauerschmidt2018dynkin}.} The integral \(\int e^{-\frac{1}{2} \Psi \Delta_W \Psi}D\Psi\) is not normalizable, we consider the graph \(\widetilde{V}\) as above, then
\begin{equation}
\label{eq-sdet}
\int \mathds{1}_{\Psi_{\delta}=0} e^{-\frac{1}{2} \widetilde{\Psi} \Delta_{\widetilde{W}} \widetilde{\Psi}}D\Psi=\int e^{-\frac{1}{2}\Psi(\Delta_W+h)\Psi} D\Psi=1.
\end{equation}

The expectation w.r.t. this density is denote \(\llbracket \cdot \rrbracket_{\widetilde{W},\Psi_{\delta}=0}\), and is called the susy free field with boundary 0 or with mass \(h\). The normalizing constant equals 1 because
\begin{equation}
\label{eq-freefield-gaussian-integral}
\int e^{-\frac{1}{2}\left< x,(\Delta_W+h)x \right>-\frac{1}{2}\left< y,(\Delta_W+h)y \right>} \prod_{i\in V}\frac{1}{2\pi}dx_idy_i=\frac{1}{\det (\Delta_W+h)}=\left( \int e^{-\left< \xi,(\Delta_W+h)\eta \right>} \prod_{i\in V}d\xi_i d\eta_i \right)^{-1}.
\end{equation}
As for the \(H^{2|2}\) field, we abusively denote \(\llbracket \cdot\rrbracket_{W}\) the non normalizable free field. The Parisi-Sourlas formula states that, if \(F( \Psi_i\Psi_j , \ i,j\in V)\) is a function of the variables \(\Psi_i\Psi_j\)s, and it decays fast enough\footnote{To simplify the exposition, we do not look for the optimal decay of \(F\) (one can check \cite{efetov1999supersymmetry} for details), in our application of Parisi--Sourlas formula, Laplace transform of local times will be enough, so we always take `suitable' test functions.}, then, denote \(A=\Delta_W+h\),
\begin{equation}
\label{eq-free-localized}
\int F e^{-\frac{1}{2} \Psi A\Psi} D\Psi=F(0).
\end{equation}
To be self-contained we include a short proof here, which we were told by M. Disertori and T. Spencer \footnote{AIM workshop on Self Interacting Processes, Supersymmetry and Bayesian Statistics 2019}.
\begin{proof}[Proof of Parisi-Sourlas formula \eqref{eq-free-localized}]
We will prove the theorem for \(|V|=1\) (general case follows in a similar manner), since \(F\) decays fast enough, we can write, for some \(\varepsilon>0\), \(F=e^{\varepsilon\Psi^2} f(\Psi^2)\), and assume that \(f\) satisfies the inverse Fourier transform formula: \(f(x)=\int \widehat{f}(k)e^{\mathrm{i}kx}dk\). Then we have
\begin{align*}
\int F e^{-\frac{1}{2} \Psi A\Psi}D\Psi&=\int \int \widehat{f}(k)e^{\mathrm{i}k\Psi^2}dk e^{\varepsilon\Psi^2-\frac{1}{2}\Psi A\Psi} D\Psi\\
&=\int \widehat{f}(k)\int e^{-\Psi (A-\mathrm{i} k-\varepsilon)\Psi} D\Psi dk\\
&=\int \widehat{f}(k) dk=f(0)=F(0),
\end{align*}
assuming \(A>\varepsilon\).
\end{proof}
Recall that if \(A\) is replaced by \(W\), the integral \(\int e^{-\frac{1}{2} \left< \Psi,\Delta_W \Psi \right>}D\Psi\) does not converge, but the above formula still holds if \(F\) is a function such that the integral converges. As an application, if we define \(z_i=\sqrt{1+\Psi_i^2 }\), and take the (fast decay) function to be
\[F=e^{\frac{1}{2} z \Delta_Wz-\left< h,z-1 \right>}\prod_{i\in V}\frac{1}{z_i}\]
for some \(h>0\), then from \(\llbracket F\rrbracket_{W}=F(0)\), we recover \(\left< 1 \right>_{\widetilde{W},\Phi_{\delta}=0}=1\), similarly, for reasonable function \(f(z)\), \(\left< f(z) \right>_{\widetilde{W},\Phi_{\delta}=0}=f(1)\). To recap, the susy integral localization is the same for free field and \(H^{2|2}\) field.

We see from the above example that, in fact, when considering susy hyperbolic sigma field and susy free field, we don't have to distinguish between \(\Phi\) and \(\Psi\), as they can be `coupled' into one (super) `probability' space. We do distinguish them in the sequel, we hope this will make the exposition more clear.

Finally, we extend our notion of susy free field on \(V\) to those with certain congruent transformation of a Markov generator, the typical case is to relate the reversible VRJP with generator \(B^u\) to \(A^u=e^{-u}B^u e^{-u}\). We can of course talk about susy free field with generator \(B^u\) (non normalizable, but can be make well defined by adding mass \(h\) as before), the expectation is denoted \(\llbracket \cdot \rrbracket_{B^u}\), now note that \(A^u= e^{-u}B^u e^{-u}\), which is not a Markov generator, we define the expectation \(\llbracket \cdot \rrbracket_{A^u}\) by
\begin{equation}
\label{eq-AuBu}
\left \llbracket F(\Psi) \right \rrbracket _{A^u}=\left \llbracket F(e^u \Psi) \right\rrbracket_{B^u}.
\end{equation}
Similarly we can add mass or pinning to these generalization of susy free fields as before. {Recall \eqref{eq-defn-Au} and \eqref{eq-defn-tilde-Au-and-tilde-Bu}.} For \(s\in \mathbb{R}\), we define
\begin{equation}
\label{eq-defn-tilde-A-s-boundary-at-delta}
\llbracket \cdot\rrbracket_{\widetilde{A}^u,\Psi_{\delta}=\theta_s0}=\int \cdot \mathds{1}_{\Psi_{\delta}=\theta_s 0}e^{-\frac{1}{2} \widetilde{\Psi} \widetilde{A}^u \widetilde{\Psi}}D\Psi,
\end{equation}
{where \(\theta_s0=(\sinh s,0,0,0)\), \(\widetilde{\Psi} \widetilde{A}^u \widetilde{\Psi}=\sum_{i,j\in\widetilde{V}}A^{u}_{ij}\Psi_i\Psi_j=\sum_{i,j\in\widetilde{V}}A^{u}_{ij}(x_ix_j+y_iy_j+\xi_i\eta_j+\xi_j\eta_i)\) and \(D\Psi=\prod_{i\in V}\frac{1}{2\pi}dx_idy_id\xi_id\eta_i\).} And for \(a\in V\), \(s\in \mathbb{R}\),
\begin{equation}
\label{eq-defn-A-s-boundary-at-a}
\llbracket \cdot \rrbracket_{A^u,\Psi_a=\theta_s 0}=\int \cdot \mathds{1}_{\Psi_a=\theta_s0} e^{-\frac{1}{2} \Psi A^u \Psi} \prod_{i\in V\setminus\{a\}}D\Psi_i,
\end{equation}
{where \(\theta_s0=(\sinh s,0,0,0)\) and \(D\Psi_i=\frac{1}{2\pi}dx_idy_id\xi_id\eta_i\).}
\section{Results}
\label{sec:org337dfe4}
The first part of our results is to provide a Bayes formula for VRJP with random initial local times, which relates the susy hyperbolic field to a susy free field. As an application we provide alternative proofs to the three isomorphism theorems and show their relations to classical isomorphism theorems.

{Since our definitions and notations can be heavy, to make the forthcoming results more reader friendly, we would like to recall that the definitions of symbols \(\nu_{i_0}^{W,z}\), \(\nu_{i_0}^{\widetilde{W},z}\), \(\left<\cdot\right>_{W,\Phi_{a}=\theta_s 0}\), \(\left<\cdot\right>_{\widetilde{W},\Phi_{\delta}=\theta_s 0}\), \(\llbracket \cdot \rrbracket_{A^u,\Psi_a=\theta_s 0}\) and \(\llbracket \cdot\rrbracket_{\widetilde{A}^u,\Psi_{\delta}=\theta_s 0}\) can be found respectively in \eqref{eq-nu-W-z}, \eqref{eq-nutilde-W-z}, \eqref{eq-defn-W-with-s-boundary}, \eqref{eq-defn-tilde-W-with-s-boundary}, \eqref{eq-defn-A-s-boundary-at-a} and \eqref{eq-defn-tilde-A-s-boundary-at-delta}.}

\begin{Theorem}[Susy Bayes formulae]
\label{orgc9bde69}
For sufficiently decaying function \(g\) such that the expectations exist, we have, for any \(a\in \widetilde{V}\) and \(b\in V\), any \(s\in \mathbb{R}\), the following hold,
\begin{equation}
\label{eq-bayes-1}
\left< \frac{z_a}{z_{\delta}} \int  g \nu_a^{\widetilde{W},z}(du) \right>_{\widetilde{W},\Phi_{\delta}=\theta_s 0}=\int \llbracket g\rrbracket_{\widetilde{A}^u,\Psi_{\delta}=\theta_s 0}\nu_a^{\widetilde{W},1}(du),
\end{equation}
\begin{equation}
\label{eq-bayes-2}
\left< \int g\ \nu_b^{W,z}(du) \right>_{W,\Phi_b=\theta_s 0}=\int \llbracket g\rrbracket_{{A}^u,\Psi_{{b}}=\theta_s0} \nu_b^{W,1}(du),
\end{equation}
where \(\theta_s0=(\sinh s,0,\cosh s,{0,0})\) in the case of \(H^{2|2}\) field and \(\theta_s0=(\sinh s,0,{0,0})\) in the case of free field.
\end{Theorem}
The following corollaries are first announced in \cite{bauerschmidt2018dynkin},\cite{bauerschmidt2019geometry}, but our formulation is slightly different.
\begin{Corollary}[BFS-Dynkin isomorphism]
\label{org55b773a}
Consider VRJP on a graph \(\widetilde{\mathcal{G}}=(\widetilde{V},\widetilde{E},\widetilde{W})\) killed at \(\delta\), let \(\zeta=\inf\{t>0,\ Y_t=\delta\}\), and \(L=L(\zeta)\) the final local time, i.e. \(L_i(\zeta)=1+\int_0^{\zeta} \mathds{1}_{Y_s=i}ds\), for any \(a,b\in V\), and any suitable function \(g\),
\begin{equation}
\label{eq-bfs-dynkin}
\int_0^{\infty} \mathbb{E}_{a,1}^{\widetilde{W}}\left( g(L)\mathds{1}_{Y_t=b,t<\zeta} \right) dt=\left< x_ax_b g(z) \right>_{\widetilde{W},\Phi_{\delta}=0}.
\end{equation}
In particular, when \(W_{\delta,i}=h\) for all \(i\),
\[ \mathbb{E}_{a,1}^{\widetilde{W}}(g(L)\mathds{1}_{Y_{\zeta^-}=b})=h \left< x_ax_b g(z) \right>_{\widetilde{W},\Phi_{\delta}=0}.\]
\end{Corollary}
\begin{Corollary}[Generalized second Ray--Knight theorem]
\label{org5b98f59}
Consider VRJP on \(\mathcal{G}=(V,E,W)\) start from \(a\in V\), for any \(s\ne 0\), let \(\tau(\gamma)=\inf\{t>0:\ L_a(t)\ge \gamma\}\), then for any suitable function \(g\),
\begin{equation}
\label{eq-second-rayknight}
\mathbb{E}^W_{a,1}(g(L(\tau(\cosh s))))=\left<  g(\theta_s z) \right>_{W,{\Phi_a=0}}.
\end{equation}
\end{Corollary}
\begin{Corollary}[Eisenbaum's isomorphism]
\label{org13a6617}
Consider VRJP on \(\widetilde{G}=(\widetilde{V},\widetilde{E},\widetilde{W})\) killed at \(\delta\),  let \(\zeta=\inf\{t>0,\ Y_t=\delta\}\), and \(L=L(\zeta)\) the final local time, i.e. \(L_i(\zeta)=1+\int_0^{\zeta} \mathds{1}_{Y_s=i}ds\), for any suitable function \(g\), for any \(s\ne 0\), any \(a\in V\),
\begin{equation}
\label{eq-eisenbaum}
\left<  \frac{z_a}{z_{\delta}}\mathbb{E}^{\widetilde{W}}_{a,z}(g(L)) \right>_{\widetilde{W},{\Phi_{\delta}=\theta_s 0}}=\left< \frac{x_a}{x_{\delta}} g(z) \right>_{\widetilde{W},{\Phi_{\delta}=\theta_s 0}}.
\end{equation}
\end{Corollary}
\begin{Theorem}[Relation to flat isomorphism theorems]
\label{org899af2f}
The equalities in Corollaries \ref{org55b773a}, \ref{org5b98f59}, \ref{org13a6617} are annealed version of the corresponding theorems for the quenched VRJP Markov process.
\end{Theorem}
The second part of our results is about reinforced loop soup. We defer the definition of reinforced loop soup to Section~\ref{sec:org4d2c65d}, intuitively speaking, it is a random collection of loops \(\mathcal{L}_t\) of VRJP trajectories, in which the number of loops increase with time \(t\). See \cite{le2011markov} for an introduction to Markov loop soup.
\begin{Theorem}[Reinforced loop soup BFS--Dynkin isomorphism]
\label{org6c0a237}
Consider the reinforced loop soup defined in \eqref{eq-def-reinforced-loop-soup-field} at time \(1\), then its occupation field equals in law to \((x^2+y^2)\) in the susy hyperbolic sigma model, more precisely,
\begin{equation}
\label{eq-soup-dynkin-susyhyp}
\mathbb{E}^{\text{soup}}(g(\widehat{\mathcal{L}_1}))=\left< g(x^2+y^2) \right>_{\widetilde{W},\Phi_{\delta}=0}.
\end{equation}
\end{Theorem}

\section{Relation between flat and hyperbolic isomorphism theorems}
\label{sec:org8ed9f79}
In this section we give alternative proofs to the three isomorphism theorems, in particular, our proofs show that the susy hyperbolic isomorphisms are in fact annealed versions of the classical isomorphism theorems. The main tool we used is a Bayes formula for VRJP with random initial local times.
\begin{proof}[Proof of Bayes formula]
All the formulae are proved by straight computations. For \eqref{eq-bayes-1}, we can combine \eqref{eq-defn-tilde-W-with-s-boundary}, \eqref{eq-nutilde-W-z} and \eqref{eq-defn-tilde-A-s-boundary-at-delta}. For \eqref{eq-bayes-2}, we combine \eqref{eq-defn-W-with-s-boundary}, \eqref{eq-nu-W-z} and \eqref{eq-defn-A-s-boundary-at-a}. {We provide detailed derivation of \eqref{eq-bayes-1} in the following and left the similar computation of \eqref{eq-bayes-2} to the reader.

Let \(F(z)=\int g(z,u)\nu_a^{\widetilde{W},z}(du)\). By \eqref{eq-defn-tilde-W-with-s-boundary}, we have that
\begin{align*}
\left< \frac{z_a}{z_{\delta}} \int  g \nu_a^{\widetilde{W},z}(du) \right>_{\widetilde{W},\Phi_{\delta}=\theta_s 0}&=\left< \frac{z_a}{z_{\delta}} F(z) \right>_{\widetilde{W},\Phi_{\delta}=\theta_s 0}\\
&=\int \frac{z_a}{z_{\delta}}F(z) \mathds{1}_{\Phi_{\delta}=\theta_s 0}e^{-\frac{1}{2} \widetilde{\Phi}\Delta_{\widetilde{W}} \widetilde{\Phi}} D\mu_V(\Phi).
\end{align*}
Recall that \(\theta_s=(\sinh s, 0, \cosh s, 0, 0)\), \eqref{eq-h22-haar} and \eqref{eq-h22-energy-mass}. We have that
\begin{multline}\label{eq-proof-bayes-1}
\left< \frac{z_a}{z_{\delta}} \int  g \nu_a^{\widetilde{W},z}(du) \right>_{\widetilde{W},\Phi_{\delta}=\theta_s 0}\\
= \int \frac{z_a}{z_{\delta}}F(z) \mathds{1}_{\Phi_{\delta}=\theta_s 0}e^{\sum_{\{i,j\}\in \tilde{E}}W_{i,j}(x_ix_j+y_iy_j-z_iz_j+\xi_{i}\eta_{j}+\xi_{j}\eta_{i}+1)} \prod_{i\in V}\frac{1}{2\pi z_i}dx_i dy_id\xi_{i}d\eta_i\\
= \int F(z) \mathds{1}_{\Phi_{\delta}=\theta_s 0}e^{\sum_{\{i,j\}\in \tilde{E}}W_{i,j}(x_ix_j+y_iy_j-z_iz_j+\xi_{i}\eta_{j}+\xi_{j}\eta_{i}+1)} \prod_{i\in \widetilde{V}\setminus\{a\}}\frac{1}{z_i}\prod_{i\in V}\frac{dx_i dy_id\xi_{i}d\eta_i}{2\pi }.
\end{multline}
Next, we use \eqref{eq-nutilde-W-z} and obtain that
\begin{equation}\label{eq-proof-bayes-2}
\begin{aligned}
F(z)&=\int g(z,u)\nu_a^{\widetilde{W},z}(du)\\
&=\int g(z,u)\mathds{1}_{u_{a}=0}e^{-\frac{1}{2}\sum_{\{i,j\}\in \widetilde{E}} W_{i,j}(e^{u_i-u_j}z_j^2+e^{u_j-u_i}z_i^2-2z_iz_j)}\sqrt{D(W,u)}\prod_{i\in \widetilde{V}\setminus \{a\}}\frac{z_i e^{-u_i}}{\sqrt{2\pi}}du_i
\end{aligned}
\end{equation}
Plugging \eqref{eq-proof-bayes-2} into \eqref{eq-proof-bayes-1}, we have that
\begin{equation}\label{eq-proof-bayes-3}
\begin{aligned}
&\left< \frac{z_a}{z_{\delta}} \int  g \nu_a^{\widetilde{W},z}(du) \right>_{\widetilde{W},\Phi_{\delta}=\theta_s 0}\\
&=\int\int g(z,u)\mathds{1}_{u_{a}=0}e^{-\frac{1}{2}\sum_{\{i,j\}\in \widetilde{E}} W_{i,j}(e^{u_i-u_j}z_j^2+e^{u_j-u_i}z_i^2-2z_iz_j)}\sqrt{D(W,u)}\prod_{i\in \widetilde{V}\setminus \{a\}}\frac{z_i e^{-u_i}}{\sqrt{2\pi}}du_i\\
&\quad \mathds{1}_{\Phi_{\delta}=\theta_s 0}e^{\sum_{\{i,j\}\in \tilde{E}}W_{i,j}(x_ix_j+y_iy_j-z_iz_j+\xi_{i}\eta_{j}+\xi_{j}\eta_{i}+1)} \prod_{i\in \widetilde{V}\setminus\{a\}}\frac{1}{2\pi z_i}\prod_{i\in V}dx_i dy_id\xi_{i}d\eta_i\\
&=\int\int g(z,u)\mathds{1}_{u_{a}=0}e^{-\frac{1}{2}\sum_{\{i,j\}\in \widetilde{E}} W_{i,j}(e^{u_i-u_j}z_j^2+e^{u_j-u_i}z_i^2)}\sqrt{D(W,u)}\prod_{i\in \widetilde{V}\setminus \{a\}}\frac{e^{-u_i}}{\sqrt{2\pi}}du_i\\
&\quad \mathds{1}_{\Phi_{\delta}=\theta_s 0}e^{\sum_{\{i,j\}\in \tilde{E}}W_{i,j}(x_ix_j+y_iy_j+\xi_{i}\eta_{j}+\xi_{j}\eta_{i}+1)} \prod_{i\in V}\frac{1}{2\pi}dx_i dy_id\xi_{i}d\eta_i.
\end{aligned}
\end{equation}
Using \(z_{i}^2=1+x_i^2+y_i^2+2\xi_i\eta_i\), we calculate the exponent in the integral:
\begin{equation}\label{eq-proof-bayes-4}
\begin{aligned}
 &-\frac{1}{2}\sum_{\{i,j\}\in \widetilde{E}} W_{i,j}(e^{u_i-u_j}z_j^2+e^{u_j-u_i}z_i^2)+\sum_{\{i,j\}\in \tilde{E}}W_{i,j}(x_ix_j+y_iy_j+\xi_{i}\eta_{j}+\xi_{j}\eta_{i}+1)\\
 &=-\frac{1}{2}\sum_{\{i,j\}\in \widetilde{E}} W_{i,j}(e^{u_i-u_j}+e^{u_j-u_i}-2)\\
 &\quad + \sum_{\{i,j\}\in \tilde{E}}W_{i,j}(x_ix_j+y_iy_j+\xi_{i}\eta_{j}+\xi_{j}\eta_{i})- \frac{1}{2}\sum_{i\in\widetilde{V}}\sum_{k\in\widetilde{V}}W_{ik}e^{u_k-u_i}(x_i^2+y_i^2+2\xi_i\eta_i)\\
 &=-\frac{1}{2}\sum_{\{i,j\}\in \widetilde{E}} W_{i,j}(e^{u_i-u_j}+e^{u_j-u_i}-2)-\frac{1}{2}\sum_{i,j}A^{u}_{ij}\Psi_{i}\Psi_{j}\\
 &=-\frac{1}{2}\sum_{\{i,j\}\in \widetilde{E}} W_{i,j}(e^{u_i-u_j}+e^{u_j-u_i}-2)-\frac{1}{2}\Psi A^{u}\Psi,
\end{aligned}
\end{equation}
where \((A^{u})_{ij}=\left\{\begin{array}{ll}
 -W_{ij}, & i\neq j\in\widetilde{V},\\
 \sum_{k\in\widetilde{V}}W_{ik}e^{u_k-u_i}, & i=j\in\widetilde{V},
\end{array}\right.\) and \(\Psi_i\Psi_j=x_ix_j+y_iy_j+\xi_{i}\eta_{j}+\xi_{j}\eta_{i}\). By \eqref{eq-nutilde-W-z}, \eqref{eq-defn-tilde-A-s-boundary-at-delta}, \eqref{eq-proof-bayes-3} and \eqref{eq-proof-bayes-4}, we obtain that
\begin{equation}
\begin{aligned}
 \left< \frac{z_a}{z_{\delta}} \int  g \nu_a^{\widetilde{W},z}(du) \right>_{\widetilde{W},\Psi_{\delta}=\theta_s 0}&=\int\int g(z,u)\mathds{1}_{\Phi_{\delta}=\theta_s 0}e^{-\frac{1}{2}\Psi A^u\Psi}\\
 &\quad\prod_{i\in V}\frac{1}{2\pi}dx_idy_id\xi_id\eta_i\nu^{\widetilde{W},1}_{a}(du)\\
 &=\int \llbracket g\rrbracket_{\widetilde{A}^u,\Psi_{\delta}=\theta_s 0}\nu_a^{\widetilde{W},1}(du).
\end{aligned}
\end{equation}
}
\end{proof}
\subsection{Proof of BFS--Dynkin's isomorphism}
\label{sec:orgdbea580}
\begin{proof}[Proof of BFS-Dynkin's isomorphism]
Denote \(\widetilde{\zeta}=\inf\{t>0:\ Z_t=\delta\}\) the time when \(Z\) is killed. Therefore, let \(S=S(\widetilde{\zeta})\) be the final local time of \(Z\),
\begin{align*}
\mathbb{E}_{a,1}^{\widetilde{W}}(g(L)\mathds{1}_{Y_{\zeta^-}=b})&=\left< z_a \mathbb{E}_{a,z}^{\widetilde{W}}(g(L)\mathds{1}_{Y_{\zeta^-}=b}) \right>_{\widetilde{W},\Phi_{\delta}=0}\\
&=\left< z_a \int  E_a^u \left(g(\sqrt{S+z^2})\mathds{1}_{Z_{\widetilde{\zeta}^-}=b }\right) \nu_a^{\widetilde{W},z}(du) \right>_{\widetilde{W},\Phi_{\delta}=0}\\
&=\int \left\llbracket E_a^u \left(g(\sqrt{S+z^2})\mathds{1}_{Z_{\widetilde{\zeta}^-}=b }\right)  \right\rrbracket_{\widetilde{A}^u,\Psi_{\delta}=0} \nu_a^{\widetilde{W},1}(du).
\end{align*}
Conditionally on \(u\), the process is Markov, so there is a classical BFS--Dynkin's isomorphism (see Appendix~\ref{sec:org82ce996} for details), which is the following
\begin{theoA}[Classical susy BFS--Dynkin's isomorphism]
\label{orgdd5023d}
Let \(Z\) be a Markov jump process on \(\mathcal{G}\) with jump rate \(\frac{1}{2}W_{i,j}e^{u_j-u_i}\) and killed when it hits \(\delta\), let \(S\) be its final local times, and let \((x,y,\xi,\eta)\) be susy free field with generator \(A^u\), with \(z^2=x^2+y^2+2\xi\eta+1\),
\begin{equation}
\label{eq-classical-susy-dynkin}
\left\llbracket E_a^u\left( g(S+z^2) \mathds{1}_{Z_{\widetilde{\zeta}^-}=b} \right)\right\rrbracket_{\widetilde{A}^u,\Psi_{\delta}=0}=W_{\delta,b}\left\llbracket e^{u_{\delta}-u_a} x_ax_b g(z^2) \right\rrbracket_{\widetilde{A}^u,\Psi_{\delta}=0}.
\end{equation}
\end{theoA}
If we apply this theorem to our previous computation and use \eqref{eq-change-root} and \eqref{eq-bayes-1}, we get that
\begin{align*}
\mathbb{E}_{a,1}^{\widetilde{W}}(g(L)\mathds{1}_{Y_{\zeta^-}=b})&=W_{\delta,b}\int \left\llbracket e^{u_{\delta}-u_a}x_ax_b g(z)  \right\rrbracket_{\widetilde{A}^u,\Psi_{\delta}=0} \nu_a^{\widetilde{W},1}(du)\\
&=W_{\delta,b}\int \nu_{\delta}^{\widetilde{W},1}(du)\left\llbracket g(z)x_ax_b  \right\rrbracket_{\widetilde{A}^u,\Psi_{\delta}=0}\\
&=W_{\delta,b}\left< x_ax_bg(z)\int \nu_{\delta}^{\widetilde{W},z}(du) \right>_{\widetilde{W},\Phi_{\delta}=0}\\
&=W_{\delta,b}\left< x_ax_bg(z) \right>_{\widetilde{W},\Phi_{\delta}=0}.
\end{align*}
The particular case when \(W_{\delta,i}=h\) for all \(i\in V\) follows immediately, since \(\zeta\) becomes an independent exponential random variable of rate \(h\).
\end{proof}
\subsection{Proof of second Ray--Knight theorem}
\label{sec:org9f3c9dd}
\begin{proof}[Proof of second Ray--Knight Theorem]
This time we will use \eqref{eq-bayes-2}. Recall that \(\tau(\gamma)=\inf\{t>0:\ L_a(t)\ge \gamma\}\), define for \(Z\), \(\sigma(\gamma)=\inf\{s>0:\ S_a(s)\ge \gamma\}\). Since we work with initial local time \(z_a=1\), we have \(L_i(\tau(\cosh s))=\sqrt{S_i(\sigma(\sinh^2 s))+z_i^2}\) for all \(i\in V\). Therefore, we have
\begin{align*}
\mathbb{E}_{a,1}^W(g(L(\tau(\cosh s))))&=\left< \mathbb{E}_{a,z}(g(L(\tau(\cosh s)))) \right>_{W,\Phi_a=0}\\
&=\left<  \mathbb{E}_{a,z}\left( g\left( \sqrt{S(\sigma(\sinh^2 s))+z^2} \right) \right) \right>_{W,\Phi_a=0}\\
&=\left< \int E_a^u\left( g\left( \sqrt{ S(\sigma(\sinh^2 s))+z^2} \right) \right) \nu_a^{W,z}(du)  \right>_{W,\Phi_a=0}\\
&=\int \left\llbracket E_a^u\left( g\left( \sqrt{S(\sigma(\sinh^2 s))+z^2} \right) \right) \right\rrbracket_{A^u,\Psi_a=0} \nu_a^{W,1}(du)
\end{align*}
Now we can apply the classical susy generalized second Ray--Knight theorem for the Markov jump process under the law \(E_a^u\), which states (see Appendix~\ref{sec:org82ce996} for details):
\begin{theoA}[Classical susy generalized second Ray--Knight theorem]
\label{org7429189}
Consider the Markov jump process on \(\mathcal{G}\) with generator \(\frac{1}{2}e^{-u} A^u e^u\), and \(\sigma(\gamma)=\inf\{s>0:\ S_a(s)\ge \gamma\}\), where \(S\) is its local times. For any \(s\ne 0\), if \((x,y,\xi,\eta)\) is a susy free field , then for any test function \(g\), with \(z^2=x^2+y^2+2\xi\eta+1\),
\begin{equation}
\label{eq-classic-susy-rayknight}
\left\llbracket E_a^u \left( g\left( S(\sigma(\sinh^2 s))+z^2 \right) \right)\right\rrbracket_{A^u,\Psi_a=0}=\left\llbracket   g(z^2) \right\rrbracket_{A^u,\Psi_{a}=\theta_s 0}
\end{equation}
\end{theoA}
Plug the theorem to our previous computations, we obtain that
\begin{align*}
\mathbb{E}_{a,1}^W(g(L(\tau(\cosh s))))&=\int \left\llbracket g(z) \right\rrbracket_{A^u,\Psi_{a}=\theta_{s}0} \nu_a^{W,1}(du)\\
&=\left\langle \int g(z)\,\nu_{a}^{W,z}(du) \right\rangle_{W,\Phi_{a}=\theta_{s}0}\\
&=\langle g(z) \rangle_{W,\Phi_{a}=\theta_{s}0}\\
&=\left< g(\theta_s z) \right>_{W,\Phi_a=0}.
\end{align*}
\end{proof}
\subsection{Proof of Eisenbaum's isomorphism}
\label{sec:org853347d}
\begin{proof}[Proof of Eisenbaum's isomorphism]
This time we will use \eqref{eq-bayes-1}, denote \(S\) the final local time of the process \(Z\) (killed at \(\delta\)),
\begin{align*}
\left< \frac{z_a}{z_{\delta}} \mathbb{E}_{a,z}^{\widetilde{W}}(g(L)) \right>_{\widetilde{W},\Phi_{\delta}=\theta_s0}&=\left< \frac{z_a}{z_{\delta}} \mathbb{E}_{a,z}\left( g\left( \sqrt{S+z^2} \right) \right) \right>_{\widetilde{W},\Phi_{\delta}=\theta_s 0}\\
&=\left< \frac{z_a}{z_{\delta}} \int E_a^u\left( g\left( \sqrt{S+z^2} \right) \right) d\nu_a^{\widetilde{W},z}(u) \right>_{\widetilde{W},\Phi_{\delta}=\theta_s0}\\
&=\int \left\llbracket E_a^u\left( g\left( \sqrt{S+z^2} \right) \right) \right\rrbracket_{\widetilde{A}^u,\Psi_{\delta}=\theta_s0} \nu_a^{\widetilde{W},1}(du)
\end{align*}
Now apply the following classical susy Eisenbaum's isomorphism (see Appendix~\ref{sec:org82ce996} for details):
\begin{theoA}[Classical susy Eisenbaum's isomorphism]
\label{orgf2011a3}
Let \(S\) be final local time of Markov process with generator \(\frac{1}{2}e^{-u}\widetilde{A}^u e^u\) on \(\mathcal{G}\) killed at \(\delta\), and \((x,y,\xi,\eta)\) a susy free field with generator \(\widetilde{A}^u\), for any \(s\ne 0\), and for suitable test function \(g\), with \(z^2=x^2+y^2+2\xi\eta+1\),
\begin{equation}
\label{eq-classical-susy-eisenbaum}
\left\llbracket E_a^u\left( g(S+z^2) \right) \right\rrbracket_{\widetilde{A}^u,\Psi_{\delta}=\theta_s 0}=e^{u_{\delta}-u_a}\left\llbracket \frac{x_a}{x_{\delta}}  g(z^2) \right\rrbracket_{\widetilde{A}^u,\Psi_{\delta}=\theta_s0}.
\end{equation}
\end{theoA}
Plug this theorem into our computation, we get that
\begin{align*}
\left< \frac{z_a}{z_{\delta}} \mathbb{E}_{a,z}^{\widetilde{W}}(g(L)) \right>_{\widetilde{W},\Phi_{\delta}=\theta_s0}&=\int e^{u_{\delta}-u_a} \left\llbracket \frac{x_a}{x_{\delta}} g(z) \right\rrbracket_{\widetilde{A}^u,\Psi_{\delta}=\theta_s 0} \nu_a^{\widetilde{W},1}(du)\\
&=\int  \left\llbracket \frac{x_a}{x_{\delta}} g(z) \right\rrbracket_{\widetilde{A}^u,\Psi_{\delta}=\theta_s 0} \nu_{\delta}^{\widetilde{W},1}(du)\\
&=\left< \frac{x_a}{x_{\delta}} g(z) \int d\nu_{\delta}^{\widetilde{W},z}(u) \right>_{\widetilde{W},\Phi_{\delta}=\theta_s0}\\
&=\left< \frac{x_a}{x_{\delta}} g(z)  \right>_{\widetilde{W},\Phi_{\delta}=\theta_s0}.
\end{align*}
\end{proof}

\section{Reinforced loop soup.}
\label{sec:org4d2c65d}
As it is shown that susy hyperbolic isomorphisms are simply annealed versions of classical isomorphisms related to free field, we can hence discuss another susy hyperbolic isomorphism in the same manner, namely the BFS--Dynkin's isomorphism for reinforced loop soup. To be self-contained, we give a brief definition of this classical object (cf.~\cite{le2011markov}) here.

Sample the environment \(\{u_i,\ i\in V,\ u_{\delta}=0\}\) according to the mixing measure \(\nu_{\delta}^{\widetilde{W},1}(du)\), from now on, fix this environment, let us first construct the quenched loop soup. To do so, consider the Markov process \(Z^u\) with jump rate \(\frac{1}{2} W_{i,j}e^{u_j-u_i}\) on \(\widetilde{V}\), killed at \(\delta\), i.e. the generator is \(L^u\) where
\begin{equation}
\label{eq-generator-loopsoup}
L^u(i,j)=\begin{cases} -\frac{1}{2} W_{i,j}e^{u_j-u_i} & i\ne j;\ i,j \in V\\ \sum_{k\in \widetilde{V}\setminus\{i\}} \frac{1}{2} W_{i,k}e^{u_k-u_i} & i=j.\end{cases}
\end{equation}
The semi-group  associated to this quenched Markov process is \(P_t^u=e^{-tL^u}\), i.e. \(P_t^u(i,j)=E_i^u(\mathds{1}_{Z_t^u=j})\), where \(E_i^u\) is the expectation of \(Z^u\) starts from \(i\).

The collection of based loops on \(V\) are admissible trajectories in \(V\) with same end points, based loop measure associated to \(L^u\) is a (usually not probability) measure on based loops. More precisely, for a given path
\[\left\{Z^u_{[0,t_1)}=i=i_0, \ Z^u_{[t_1,t_2)}=i_1 ,\ldots,\ Z^u_{[t_{k-1},t_k)}=i_{k-1},\ Z^u_{[t_k,t]}=j=i_k\right\},\]
which we will denote
\[i\xrightarrow[]{t_1}i_1\xrightarrow[]{t_2-t_1}i_2\rightarrow\ldots \rightarrow i_{k-1}\xrightarrow[]{t_k-t_{k-1}}j\xrightarrow[]{t-t_k};\]
The trivial path (that makes no jumps) is denoted \(i_0\xrightarrow[]{t}\). In the remainder, the path measure is defined, for fixed \(t>0\),
\begin{multline}
\label{equation-path-measure}
P_t^{u,(i,j)}\left(i\xrightarrow[]{t_1}i_1\xrightarrow[]{t_2-t_1}i_2\rightarrow\ldots \rightarrow i_{k-1}\xrightarrow[]{t_k-t_{k-1}}j\xrightarrow[]{t-t_k}  \right)\\=L^u_{i,i_1}L^u_{i_1,i_2}\cdots L^u_{i_{k-1},j}e^{-t_1 L^u_{i,i}-(t_2-t_1)L^u_{i_1,i_1}-\cdots-(t_k-t_{k-1})L^u_{i_{k-1},i_{k-1}} -(t-t_k)L^u_{j,j}}dt
\end{multline}
where \(dt=dt_1\cdots dt_k\); and \(P_t^{u,(i_0,i_0)}({i_0\xrightarrow[]{t}})=e^{-t L_{i_0,i_0}^u}\). It is a measure on paths up to time \(t\) (and \(k\ge 0\) is arbitrary), with \(Z^u_0=i,\ Z^u_t=j\). In particular, \(\sum_j P_t^{(i,j)}\) is the marginal of the probability of the Markov process up to time \(t\).

A based loop with base point \(i\) is a path of the following form
\[i\xrightarrow[]{t_1}i_1\xrightarrow[]{t_2-t_1}i_2\rightarrow\ldots \rightarrow i_{k-1}\xrightarrow[]{t_k-t_{k-1}}i\xrightarrow[]{t-t_k}\]
The based loop measure is a measure on the collection of based loops (endowed with the sigma algebra generated by finite dimensional marginals):
\[\left\{i\xrightarrow[]{t_1}i_1\xrightarrow[]{t_2-t_1}i_2\rightarrow\ldots \rightarrow i_{k-1}\xrightarrow[]{t_k-t_{k-1}}i\xrightarrow[]{t-t_k}, \begin{array}{l} k\ge 0, \ t>0, i_0,i_1 ,\ldots,i_k\in V,\\  0<t_1<t_2<\cdots <t_k<t\end{array}\right\},\]
defined by
\begin{equation}
\mu^{b,u}=\sum_{i\in V}\int_0^{\infty} \frac{1}{t}P_t^{u,(i,i)}dt.
\end{equation}
Note that, \(P_t^{u,(i,i)}\) is a non-normalized bridge measure from \(i\) to \(i\).

Once we have a measure on based loops, we can forget the base point of a based loop, and work on the equivalence class of based loops. The loop measure \(\mu^u\) is the image measure of based loop measure w.r.t. this equivalence relation. The observables on loop soup that we will consider are loop functionals, so all the computations can be done on based loops. For example, a particularly interesting loop functional is \(\widehat{\ell}\), which is the occupation time of a loop (or based loop, which makes no difference):
\[\widehat{\ell}_i=\int_0^t \mathds{1}_{\ell(s)=i}ds,\ \forall i\in V.\]
In the sequel, we will not distinguish based loops and loops.

The loop soup \((\mathcal{L}_{\alpha}^u,\ \alpha>0)\) is a Poisson point process on the space of loops of intensity \(\mu^u\). In particular, for any measurable subset \(F\) of loops, the number of loops of \(\mathcal{L}_{\alpha}^u\) in \(F\) is Poisson distributed with parameter \(\alpha\mu^u(F)\). The occupation field \(\widehat{\mathcal{L}_{\alpha}^u}\) is defined by \(\widehat{\mathcal{L}_{\alpha}^u}=\sum_{\ell\in \mathcal{L}_{\alpha}^{u}}\widehat{\ell}\). The distribution of the loop soup \(\mathcal{L}^u\) is denoted \(\mathbb{E}^{\text{soup},u}\).

The following theorem for loop soup is well known (e.g. p48. Remark (c) after Theorem~2 in \cite{le2011markov}):
\begin{theoA}[BFS--Dynkin's isomorphism for loop soup]
\label{org1872e7d}
With the definition in this section, the occupation field of \(\mathcal{L}_1^u\) has the same distribution of the sum of squares of two Gaussian free fields with generator \(A^u\). That is, for any suitable test function \(g\),
\begin{equation}
\label{eq-loopsoup-dynkin-classic}
\mathbb{E}^{\text{soup},u}(g(\widehat{\mathcal{L}_1^u}))=\int g(x^2+y^2) e^{-\frac{1}{2} (x A^u x +yA^u y)} \det (A^u) dx dy.
\end{equation}
\end{theoA}
With the trick in appendix, we can translate this theorem into its susy free field version.
\begin{Corollary}[BFS--Dynkin's isomorphism for loop soup, susy free field version]
\label{org9c310d8}
For suitable test function \(g\), with \(z^2=x^2+y^2+2\xi\eta+1\),
\begin{equation}
\label{eq-dynkin-soup-susy-free}
\left\llbracket \mathbb{E}^{\text{soup},u} (g(\widehat{\mathcal{L}_1^u}+2\xi\eta+1))\right\rrbracket_{\widetilde{A}^u,\Psi_{\delta}=0}=\left\llbracket g(z^2) \right\rrbracket_{\widetilde{A}^u,\Psi_{\delta}=0}.
\end{equation}
\end{Corollary}
Motivated by the above theorem, we define occupation field for reinforced loop soup to be the annealed version of \(\widehat{\mathcal{L}_{\alpha}^u}\) under \(\nu_{\delta}^{\widetilde{W},1}(du)\). That is
\begin{equation}
\label{eq-def-reinforced-loop-soup-field}
\mathbb{E}^{\text{soup}}(g(\widehat{\mathcal{L}_{\alpha}}))=\int \mathbb{E}^{\text{soup},u}(g(\widehat{\mathcal{L}_{\alpha}^u})) d\nu_{\delta}^{\widetilde{W},1}(u).
\end{equation}
\begin{Remark}[Time change in reinforced loop soup]
\label{orge9d43cd}
One would ask for a natural definition of reinforced loop soup, in the manner of throwing loops on the graph, each loop will update the occupation field when it is threw, thus results in changing the initial local time for the next loop. A definition with this kind of feature can be carried out by using Remark 21 in~\cite{le2011markov}, by cutting a Markov paths into loops. We plan to develop these aspects in a further work.
\end{Remark}

\begin{proof}[Proof of susy hyperbolic loop soup BFS--Dynkin's isomorphism]
\begin{align*}
\mathbb{E}^{\text{soup}}(g(\widehat{\mathcal{L}_1}))&=\int \mathbb{E}^{\text{soup},u}( g(\widehat{\mathcal{L}_1^u}))d\nu_{\delta}^{\widetilde{W},1}(u)\\
&=\int \left\llbracket g(x^2+y^2) \right\rrbracket_{\widetilde{A}^u,\Psi_{\delta}=0} d\nu_{\delta}^{\widetilde{W},1}(u)\\
&=\left< g(x^2+y^2) \int d\nu_{\delta}^{\widetilde{W},z} \right>_{\widetilde{W},\Phi_{\delta}=0}=\left< g(x^2+y^2) \right>_{\widetilde{W},\Phi_{\delta}=0}.
\end{align*}
\end{proof}

\appendix
\section{the susy free field versions of classical isomorphism theorems}
\label{sec:org82ce996}
The standard BFS--Dynkin's isomorphism, generalized second Ray--Knight theorem and Eisenbaum's isomorphism for reversible Markov processes are related to Gaussian free field, see \cite{sznitman2012topics} for a thorough discussion. Here we would like to explain their slightly more general susy versions. Basically, we use the following two integrals: for a symmetric real M-matrix \(A>0\) of size \(V\times V\), say for some \(h>0\),
\[A_{i,j}=\begin{cases} -W_{i,j} & i\ne j\\ W_{i,i}+h & i=j\end{cases}.\]
We have
\begin{equation}
\label{eq-realgff}
\int_{\mathbb{R}^V} e^{-\frac{1}{2} yAy}\prod_{i\in V} \frac{dy_i}{\sqrt{2\pi}}=\frac{1}{\sqrt{\det{A}}},
\end{equation}
and
\begin{equation}
\label{eq-fermiongff}
\int e^{-\xi A \eta}d\xi d\eta=\det A.
\end{equation}
Now, consider for example the BFS--Dynkin's isomorphism, which states that, for a Markov process \(X_t\) of generator \(A\), with law \(E_a^A\), killed at rate \(h\), at time \(\zeta\). Let \(L\) be its final local times, and \(x\) is a real Gaussian free field of generator \(A\), then, we have the equality in Laplace transform:
\begin{align*} \int \sqrt{\det A} E_{a,b}^A \left(e^{-kL-\frac{1}{2} k x^2}\right) e^{-\frac{1}{2} xAx} \prod_{i\in V}\frac{dx_i}{\sqrt{2\pi}}& =(A+k)^{-1}(a,b) \sqrt{\frac{\det A}{\det (A+k)}}\\
&=\int \sqrt{\det A} x_ax_b e^{-\frac{1}{2} k x^2} e^{-\frac{1}{2} xAx}\prod_{i\in V}\frac{dx_i}{\sqrt{2\pi}}.\end{align*}
where \(E_{a,b}^A(\cdot)=\frac{1}{h}E_a^A(\cdot \mathds{1}_{X_{\zeta^-}}=b)\). Now on both the LHS and RHS of the above equality, we can replace the \(\sqrt{\det A}\) by
\[\sqrt{\det A}=\frac{1}{\sqrt{\det A}} \cdot \det A\]
and use \eqref{eq-realgff} and \eqref{eq-fermiongff} to write
\[\sqrt{\det A}=\int e^{-\frac{1}{2} yAy-\xi A\eta} \prod_{i\in V}\frac{dy_i}{\sqrt{2\pi}}d\xi_id\eta_i,\]
Now if we replace \(A\) by \(e^{-u}A^ue^u\) as in \eqref{eq-AuBu}, we recover Theorem \ref{orgdd5023d}. Similar arguments provide susy counterparts of the other two theorems.

\section{Another proof of BFS--Dynkin's isomorphism à la Feynman-Kac}
\label{sec:orge6ba3f6}
In this appendix we provide another proof of BFS--Dynkin's isomorphism, i.e. Theorem \ref{org55b773a}. This proof is done in a more classic way, i.e. it uses some kind of Feynman--Kac idea. We use the notion of trajectory density introduced in Lemma 1 of \cite{zeng2013vertex}, by this lemma the trajectory density of a VRJP \((Y_t)\) on \(\mathcal{G}\) with starting point \(a\in V\) up till some \(t>0\) equals
\[d_{\sigma}=\mathds{1}_{0<t_1<\cdots <t_n<t}\prod_{k=1}^n W_{i_{k-1},i_k}L_{i_k}(t_{k-1}) e^{-\int_0^t \sum_{j\in V}W_{Y_u,j}L_j(u)du}dt_1\cdots dt_n\]
where \(\sigma\) is the trajectory up to time \(t\), defined by
\[\sigma=\{Y_{[0,t_1)}=i_0=a,\ Y_{[t_1,t_2)}=i_1 ,\ldots,Y_{[t_{n-1},t_n)}=i_{n-1}, \ Y_{[t_n,t]=i_n=b}\}\]
As initial local times of \(Y\) equal 1 everywhere,
\[k\cdot \ell(\sigma):=k\ell(\sigma)=\int_0^t k_{Y_u}du.\]
By Lemma 1 of \cite{zeng2013vertex}, we have
\begin{align*}
\int_0^{\infty} \mathbb{P}_{a,1}^W\left(e^{-k \ell(t)} \mathds{1}_{Y_t=b}\right) he^{-ht}dt&=\int_0^{\infty} \sum_{n\ge 0}\sum_{i_1 ,\ldots,i_{n-1}} \int \mathds{1}_{i_0=a} e^{k\ell(\sigma)} \mathds{1}_{i_n=b} d_{\sigma}   h e^{-ht}dt\\
&=h\int_0^{\infty} \sum_{n\ge 0}\sum_{i_1 ,\ldots,i_{n-1}} \int \mathds{1}_{i_0=a,\ i_n=b} \widetilde{d_{\sigma}} dt
\end{align*}
where
\[\widetilde{d_{\sigma}}=\mathds{1}_{0<t_1<\cdots <t_n<t}\prod_{k=1}^n W_{i_{k-1},i_k}L_{i_k}(t_{k-1}) e^{-\int_0^t \sum_{j \in V\cup\{\iota\}}W_{Y_u,j}L_j(u)du}dt_1\cdots dt_n.\]
Here, \(\iota\) is a cemetery point and we set \(W_{i,\iota}=k_i\) for \(i\in V\). In the above way, the term \(e^{-k\ell(t)}\) is absorbed into the trajectory density. On the other hand, let \(\delta\) be another cemetery point such that \(W_{i,\delta}=h\) for \(i\in V\cup\{\iota\}\). Then the hitting time \(T_{\delta}\) of the vertex \(\delta\) has probability density \(he^{-ht}\). Let us denote by \(\mathbb{P}_{a,1}^{W,k}\) the law of the VRJP with one cemetery point \(\iota\) and by \(\mathbb{P}_{a,1}^{W,k,h}\) the law of the VRJP with two cemetery point \(\delta\) and \(\iota\).
\begin{align*}
 \int_0^{\infty} \mathbb{P}_{a,1}^W\left(e^{-k \ell(t)} \mathds{1}_{Y_t=b}\right) he^{-ht}dt=&\int_0^{\infty} \mathbb{P}_{a,1}^{W,k}\left(Y_t=b,t<T_{\iota}\right) he^{-ht}dt\\
 =&\mathbb{P}_{a,1}^{W,k,h}(Y_{T_{\delta}^-=b},\ T_{\delta}<T_{\iota}),
\end{align*}
where \(T_{\iota}=\inf\{t>0:\ Y_t=\iota\}\). Now this probability is invariant under time change as it only concerns the discrete time chain (a.k.a. the skeleton) associated to the continuous time process, and by Theorem 2.(ii) in \cite{ST15} we know that the skeleton is a mixture of Markov chain with conductance \(W_{i,j}e^{u_i+u_j}\) where \(u\) is sampled according to \(\nu^{\widetilde{\widetilde{W}},1}_{a}(du)\), where \(\widetilde{\widetilde{W}}\) is the extended conductance \(W\) on \(V\cup\{\iota,\delta\}\).

For a Markov jump process \(Z\) on \(V\cup\{\iota,\delta\}\) with a generator \(B\), by considering the corresponding discrete Markov chain, we have that
\[\mathbb{P}_{a}(Z_{T_{\delta}^-=b},\ T_{\delta}<T_{\iota})=(B|_{V\times V})^{-1}(a,b)B(b,\delta).\]
In our case, \(B\) is the generator \(B^u\) related to the conductance \(W_{i,j}e^{u_i+u_j}\) for \(i,j\in V\cup\{\iota,\delta\}\) and \(B(b,\delta)=he^{u_{\delta}+u_b}\). Hence, we have that
\[\mathbb{P}_{a,1}^{W,k,h}(Y_{T_{\delta}^-=b},\ T_{\delta}<T_{\iota})=h\int (B^u|_{V\times V})^{-1}(a,b)e^{u_a+u_b}e^{u_{\delta}-u_{a}} \nu^{\widetilde{\widetilde{W}},1}_{a}(du).\]
Performing the change of variables \(v_i=u_i-u_{\delta}\), we have that
\[\mathbb{P}_{a,1}^{W,k,h}(Y_{T_{\delta}^-=b},\ T_{\delta}<T_{\iota})=h\int (B^v|_{V\times V})^{-1}(a,b)e^{v_a+v_b} \nu^{\widetilde{\widetilde{W}},1}_{\delta}(dv).\]
By the explicit Laplace transform of \(\left( 2\beta_i:=\sum_{k\in V\cup\{\iota,\delta\}}W_{i,k}e^{u_k-u_i} \right)_{i\in V}\), e.g (5.4) of \cite{SZ15}, we see that the marginal of \((\beta_i)_{i\in V}\) are the same under \(\nu_{\delta}^{\widetilde{\widetilde{W}},1}\) and \(\nu_{\delta}^{\widetilde{W}^{k+h},1}\), it is not hard to check by a change of variable that
\[\int (B^v|_{V\times V})^{-1}(a,b)e^{v_a+v_b} \nu^{\widetilde{\widetilde{W}},1}_{\delta}(dv)=\int (B^u|_{V\times V})^{-1}(a,b)e^{u_a+u_b} \nu^{\widetilde{W}^{k+h},1}_{\delta}(du),\]
where \(\widetilde{W}^{k+h}_{i,j}=W_{i,j}\) for \(i,j\in V\) and \(\widetilde{W}^{k+h}_{i,\delta}=k_i+h\).
Now, in the horospherical coordinate (See e.g. appendix of \cite{bauerschmidt2018dynkin}) of \(H^{2|2}\) model, for the effective bosonic measure on the variables \((u,s)\), \(y=e^us\), conditionally on \(u\), \(s\) is a GFF with generator \(W_{i,j}e^{u_j+u_i}\). Besides, the distribution of \(u\) under \(\left<\cdot\right>_{\widetilde{W}^{k+h},\Phi_{\delta}=0}\) is \(\nu^{\widetilde{W}^{k+h},1}_{\delta}(du)\). Hence, by Wick's formula,
\[\left< y_ay_b \right>_{\widetilde{W}^{k+h},\Phi_{\delta}=0}=\int (B^u|_{V\times V})^{-1}(a,b)e^{u_a+u_b} \nu^{\widetilde{W}^{k+h},1}_{\delta}(du).\]
Finally, note that
\[\left< y_ay_b \right>_{\widetilde{W}^{k+h},\Phi_{\delta}=0}=\left< y_ay_be^{-\left< k,z-1 \right>} \right>_{\widetilde{W},\Phi_{\delta}=0}=\left< x_ax_be^{-\left< k,z-1 \right>} \right>_{\widetilde{W},\Phi_{\delta}=0},\]
therefore,
\[\int_0^{\infty} \mathbb{P}_{a,1}^W\left(e^{-k \ell(t)} \mathds{1}_{Y_t=b}\right) he^{-ht}dt=h\left< x_ax_b e^{-\left< k,z-1 \right>} \right>_{\widetilde{W},\Phi_{\delta}=0}.\]

\bibliographystyle{plain}
\bibliography{refs}
\end{document}